\documentclass[12pt]{article}
\usepackage{amsmath}
\usepackage{amssymb}
\usepackage[dvips]{graphicx}
\font\eightmsb=msbm10 scaled 800

\def\eps{\varepsilon}
\font\tencmmib=cmmib10 \skewchar\tencmmib '60
\newfam\cmmibfam
\textfont\cmmibfam=\tencmmib

\font\eightmsb=msbm10 scaled 1100

\def\Bbbb#1{\hbox{\eightmsb#1}}
\def\bbox{\quad\hbox{\vrule \vbox{\hrule \vskip2pt \hbox{\hskip2pt
\vbox{\hsize=1pt}\hskip2pt} \vskip2pt\hrule}\vrule}}
\def\lessim{\ \lower4pt\hbox{$
\buildrel{\displaystyle <}\over\sim$}\ }
\def\gessim{\ \lower4pt\hbox{$\buildrel{\displaystyle >}
\over\sim$}\ }

\def\P{{\rm P}}

\def\eps{{\varepsilon}}

\def\qed{\hfill\break\rightline{$\bbox$}}
\newcommand{\e}{{\Bbbb E}}
\newcommand{\p}{{\Bbbb P}}
\newcommand{\Reals}{\mathbb{R}}

\newcommand{\F}{{\cal F}}

\newtheorem{proposition}{Proposition}
\newtheorem{lemma}{Lemma}
\newtheorem{theorem}{Theorem}
\newtheorem{corollary}{Corollary}
\makeatletter
\@addtoreset{equation}{section}

\makeatother

\begin{document}

\title{
Symmetrization approach to concentration inequalities
for empirical processes.}

\author{
Dmitry Panchenko
\\
{\it Massachusetts Institute of Technology}\\
}
\date{}

\maketitle

\begin{abstract}
We introduce a symmetrization technique that allows us to
translate a problem of controlling the deviation of 
some functionals on a product space from their mean
into a problem of controlling the deviation between two
independent copies of the functional. 
As an application we give a new easy proof of Talagrand's
concentration inequality for empirical processes, where
besides symmetrization we use only Talagrand's concentration
inequality on the discrete cube $\{-1,+1\}^n.$ 
As another application of this technique we prove
new Vapnik-Chervonenkis type inequalities.
For example, for VC-classes of functions we prove
a classical inequality of Vapnik and Chervonenkis 
only with normalization by the sum of variance and sample variance. 
\end{abstract}

\vskip 5mm

\hfill\break
\hfill\break
\hfill\break
\vfill\break

\section{Introduction and main results.}

Let us consider a measurable space $\Omega$ with probability
measure $\mu,$ and the corresponding product space 
$(\Omega^n,\mu^n).$ Given a class of measurable functions
$\F=\{f:\Omega\to \Reals\},$ we consider a functional
$$
Z(x)=\sup_{\F}\sum_{i=1}^{n} f(x_i)
$$
where $x=(x_1,\ldots,x_n)\in\Omega^n,$ which is usually 
called an empirical process. To avoid measurability 
problems we will assume that $\F$ is countable, or even finite.
Our main interest is to study the deviation inequalities 
for this (or similar) functional from its mean.
The main observation of this paper is that this problem
can be translated into a problem of studying 
$Z(x)-Z(y),$ where $y$ lives on a separate copy of $\Omega^n.$
This new problem turns out to be easier, at least in the examples
we have in mind here, as it can be handled with Talagrand's 
convex distance inequality on $\{-1,+1\}^n$
which is the simplest case of 
convex distance inequality (see Talagrand (1995)).

As a first example of application of this technique we will
give an easy proof of Talagrand's concentration inequality
for $Z(x).$ As a second example, we will prove new Vapnik-Chervonenkis
type inequalities.

Let us start by proving the main result that will allow us
to implement the mentioned symmetrization.
For $x\in\Reals$ 
we will denote $(x)_{+}=\max(x,0).$

\begin{lemma}
If $\xi$ and $\nu$ are r.v.s such that for any number
$a\in\Reals$ and a function $\phi(x)=(x-a)_{+}$
$$
\e\phi(\xi)\leq \e\phi(\nu)
$$
and for some $\Gamma\geq 1, \gamma>0$ and for all $t\geq 0$
$$
\p(\nu \geq t)\leq\Gamma e^{-\gamma t},
$$
then for all $t\geq 0$
$$
\p(\xi \geq t)\leq\Gamma e^{1-\gamma t}.
$$
\end{lemma}
{\bf Proof.}
Let $\phi(x)=(x-a)_{+}$ for some $a\in\Reals$ that
will be chosen later.
Note that $\phi$ is nondecreasing. 
For $t>0$ we can write
\begin{eqnarray*}
&&
\p(\xi\geq t)\leq
\frac{\e\phi(\xi)}{\phi(t)}\leq
\frac{\e\phi(\nu)}{\phi(t)}
=
\frac{1}{\phi(t)}\Bigl(
\phi(0)+\int_{0}^{\infty}\phi'(x)\p(\nu\geq x)dx
\Bigr)
\\
&&
\leq
\frac{1}{\phi(t)}\Bigl(
\phi(0)+\Gamma\int_{0}^{\infty}\phi'(x)e^{-\gamma x}dx
\Bigr),
\end{eqnarray*}
where we used integration by parts.
Since $\Gamma\geq 1,$ we can assume that $t\geq \gamma^{-1}.$
Take
$$
a=t-\frac{1}{\gamma},\,\,\,\,\,
\phi(x)=\Bigl(x-t+\frac{1}{\gamma}\Bigr)_{+}.
$$
Then $\phi(t)=\gamma^{-1},$ $\phi(0)=0$ and
$$
\int_{0}^{\infty}\phi'(x)e^{-\gamma x}dx=
\int_{t-\gamma^{-1}}^{\infty}e^{-\gamma x}dx=
\gamma^{-1}e^{1-\gamma t},
$$
which gives 
$\p(\xi \geq t)\leq\Gamma e^{1-\gamma t}.$
\qed

It is clear that the Lemma can be stated in more generality,
for instance, we could consider the case of tails 
$\Gamma e^{-\gamma t^{\alpha}}$ for $\alpha>0.$ 
But it is irrelevant for the applications of this paper.
The main consequence is given by the following corollary.

\begin{corollary}
Let $\xi_i(x,y):\Omega^n\times\Omega^n\to R,$ 
$1\leq i\leq 3$ be measurable functions defined on two copies
of $\Omega^n$ and let
$$
\xi_i^{\prime}(x) = \int_{\Omega^n} \xi_{i}(x,y) d\mu^n(y).
$$
If $\xi_3\geq 0$ and for all $t\geq 0$
$$
\mu^{2n}(\xi_{1}\geq \xi_2 + (\xi_3 t)^{1/2} )\leq \Gamma e^{-\gamma t},
$$
then for all $t\geq 0$
$$
\mu^{n}(\xi_{1}^{\prime}\geq 
\xi_{2}^{\prime} + (\xi_{3}^{\prime} t)^{1/2} )\leq \Gamma e^{1-\gamma t}.
$$
\end{corollary}

{\bf Proof.}
Since $\sqrt{ab}=\inf_{\delta>0}(\delta a + b/(4\delta))$
we can rewrite the events
$$
\Bigl\{\xi_{1}\geq \xi_2 + (\xi_3 t)^{1/2}\Bigr\}=
\Bigl\{\sup_{\delta>0}4\delta(\xi_{1} - \xi_2 -\delta \xi_3)\geq t 
\Bigr\}
$$
and, similarly,
$$
\Bigl\{\xi_{1}^{\prime}\geq \xi_{2}^{\prime} + (\xi_{3}^{\prime} t)^{1/2}
\Bigr\}
=
\Bigl\{\sup_{\delta>0}4\delta(\xi_{1}^{\prime} - \xi_{2}^{\prime}
 -\delta \xi_{3}^{\prime})\geq t \Bigr\}. 
$$
Let us denote
$$
\xi=\sup_{\delta>0}4\delta(\xi_{1} - \xi_2 -\delta \xi_3),\,\,\,
\nu=\sup_{\delta>0}4\delta(\xi_{1}^{\prime} - \xi_{2}^{\prime}
 -\delta \xi_{3}^{\prime}).
$$
Clearly,
$$
\nu=
\sup_{\delta>0}
\int 4\delta(\xi_{1} - \xi_2 -\delta \xi_3) d\mu^n(y)\leq
\int \xi d\mu^n(y),
$$
and, thus, by Jensen's inequality,
for any nondecreasing convex funcion $\phi$
$$
\int \phi(\nu) d\mu^{n}(x) \leq
\int \phi\Bigl(\int \xi d\mu^n(y)\Bigr) d\mu^{n}(x)\leq
\int \phi(\xi) d\mu^{n}(x)d\mu^n(y).
$$
Lemma 1 implies the result.
\qed

As we mentioned above, besides the symmetrization of Corollary 1
we will need Talagrand's convex distance inequality, which we
will formulate now.

Consider the space $\{0,1\}^n$ with uniform measure $\p_{\eps}.$
If $\eps\in\{0,1\}^n$ and 
${\cal A}\subseteq \{0,1\}^n,$
denote
$$
U_{\cal A}(\eps)=\{(s_i)_{i\leq n}\in\{0,1\}^n , 
\exists \eps'\in{\cal A} , s_i=0 \Rightarrow \eps_{i}'=\eps_{i}\}.
$$
Denote the "convex hull" distance between the point 
$\eps$ and a set $\cal A$ as
$$
f_c({\cal A},\eps)=\inf\{|s| : s\in\mbox{conv}U_{\cal A}(\eps)\},
$$
where $|s|$ denotes the Euclidean norm of $s.$
The concentration inequality of Talagrand (Theorem 4.3.1 in \cite{Ta1})
states the following.

\begin{proposition}
For any $\alpha\geq 0$
\begin{equation}
\p_{\eps}(f_{c}^{2}({\cal A},\eps)\geq t)\leq
\frac{1}{\p_{\eps}({\cal A})^{\alpha}}\exp\Bigl\{
-\frac{\alpha}{\alpha+1}t\Bigr\}.
\label{T1}
\end{equation}
\end{proposition}

{\bf Remark.} In \cite{Ta1} this result was formulated
for $\alpha\geq 1,$ but it was proven (and used) for 
$\alpha\geq 0.$

The main feature of this distance is that if
$f_{c}^{2}({\cal A},\eps)\leq t,$ then (Theorem 4.1.2 in \cite{Ta1})
\begin{equation}
\forall (\lambda_i)_{i\leq n}\,\,\,\,\, 
\exists \eps'\in{\cal A}\,\,\,\,\,\,\,\,
\sum_{i=1}^{n}\lambda_i I(\eps_{i}'\not = \eps_i)\leq
(t\sum_{i=1}^{n}\lambda_i^2)^{1/2}.
\label{T2}
\end{equation}
We will start by giving a new proof of Talagrand's 
concentration inequality for empirical processes.
\vspace{0.5cm}

\section{Talagrand's concentration inequality 
for empirical processes.}

For simplicity of notations from now on we will write $\p$ to
denote any probability measure, and $\p_{\xi}$ to specify the
distribution on the space of random variable $\xi,$ with all
other variables fixed. Similarly, to denote the expectation
we will write $\e$ and $\e_{\xi}.$

Let us define a mixed uniform variance as
\begin{equation}
V=\e_{y}\sup_{f\in\F}\sum_{i=1}^{n}(f(x_i)-f(y_i))^2.
\label{univar}
\end{equation}
In a sense, $V$ is a uniform version of the sum of variance
and sample variance, since in the case when $\F$ consists of
one function, this is exactly what it is. Clearly, $V$ is a 
function of $x.$
The following theorem holds.

\begin{theorem}
Let $V$ be defined by (\ref{univar}). Then
for any $\alpha > 0$
$$
\p\Bigl(\sup_{f\in \F}\sum_{i=1}^n f(x_i)\geq 
\e\sup_{f\in\F}\sum_{i=1}^n f(x_i)
+2\sqrt{Vt}\Bigr)\leq 2^{\alpha+1}
\exp\Bigl\{1-\frac{\alpha}{\alpha+1}t\Bigr\}
$$
and
$$
\p\Bigl(\sup_{f\in \F}\sum_{i=1}^n f(x_i)\leq 
\e\sup_{f\in\F}\sum_{i=1}^n f(x_i)
-2\sqrt{Vt}\Bigr)\leq 2^{\alpha+1}
\exp\Bigl\{1-\frac{\alpha}{\alpha+1}t\Bigr\}
$$
\end{theorem}

{\bf Remark.} One can optimize the bound over $\alpha,$
which would give that for $t\geq\log 2,$ the bound 
can be written as 
$2\exp\{1-(\sqrt{t}-\sqrt{\log 2})^2\}.$

{\bf Proof.} 
We will only prove the upper tail, since the proof of the lower
tail is exactly the same, once one switches $Z$ and $EZ.$
Since
$$
\e\sup_{f\in\F}\sum_{i=1}^n f(x_i)=
\e_{y}\sup_{f\in\F}\sum_{i=1}^n f(y_i)
$$
Corollary 1 implies that it is enough to prove that
$$
\p\Bigl(\sup_{f\in \F}\sum_{i=1}^n f(x_i)\geq 
\sup_{f\in \F}\sum_{i=1}^n f(y_i)
+2\sqrt{Wt}\Bigr)\leq 2^{\alpha+1}
\exp\Bigl\{-\frac{\alpha}{\alpha+1}t\Bigr\},
$$
where $W =\sup_{f\in\F}\sum_{i=1}^{n}(f(x_i)-f(y_i))^2.$
For any $(x_1,\ldots,x_n,y_1,\ldots,y_n),$ let $\Pi$ be
the set of permutations of these coordinates such that, 
for each $1\leq i\leq n,$ $\pi(x_i),\pi(y_i)\in\{x_i,y_i\},$
and let $\p_{\pi}$ denote the uniform probability measure
on $\Pi.$ Since the above probability is invariant with respect
to any $\pi\in\Pi,$
it is enough to show that for any fixed 
$x=(x_1,\ldots,x_n)$ and $y=(y_1,\ldots,y_n)$ the probability
over permutations
$$
\p_{\pi}\Bigl(\sup_{f\in \F}\sum_{i=1}^n f(z_i^1)
\geq \sup_{f\in \F}\sum_{i=1}^n f(z_i^2)
+2\sqrt{Wt}\Bigr)\leq 2^{\alpha+1}
\exp\Bigl\{-\frac{\alpha}{\alpha+1}t\Bigr\},
$$
where $z_i^1=\pi(x_i)$ and $z_i^2=\pi(y_i).$ 
Note that $W$ is invariant under permutations.
We can rewrite it differently in terms of an i.i.d.
Bernoulli sequence $\eps=(\eps_1\ldots,\eps_n),$
i.e. $\p(\eps_i=0)=\p(\eps_i=1)=1/2.$
Namely, we can write
$$
f(z_i^1)=f(y_i)+\eps_i(f(x_i)-f(y_i)),\,\,\,
f(z_i^2)=f(x_i)-\eps_i(f(x_i)-f(y_i)),
$$
and instead of permutations look at the distribution $\p_{\eps}$ 
of $\eps.$
For any $f\in \F$ let us denote $c_f=\sum f(y_i),$
$c_f^{\prime}=\sum f(x_i),$ and $f_i=(f(x_i)-f(y_i)).$
Then, we need to prove that
$$
\p_{\eps}\biggl(
\sup_{f\in \F} \Bigl(c_f+\sum_{i=1}^n \eps_i f_i\Bigr) \geq 
\sup_{f\in \F} \Bigl(c_f^{\prime}-\sum_{i=1}^n \eps_i f_i\Bigr) +
2\Bigl(t\ \sup_{f\in \F}\sum_{i=1}^n f_i^2\Bigr)^{1/2}  
\biggr)\leq 2^{\alpha+1}
\exp\Bigl\{-\frac{\alpha}{\alpha+1}t\Bigr\}.
$$
But this is an easy consequence of Proposition 1.
Let us consider the functionals
$$
\Phi(\eps)=\sup_{f\in \F} \Bigl(c_f+\sum_{i=1}^n \eps_i f_i\Bigr),\,\,\,
\Phi^{\prime}(\eps)=\sup_{f\in \F} 
\Bigl(c_f^{\prime}-\sum_{i=1}^n \eps_i f_i\Bigr).
$$
They are both convex, with the Lipschitz norm bounded by
$$
\|\Phi\|_{L}, \|\Phi^{\prime}\|_{L}\leq
\Bigl(\sup_{f\in \F}\sum f_i^2\Bigr)^{1/2}.
$$
Also, by symmetry, they have the same median,
$M=M(\Phi)=M(\Phi^{\prime})$ with respect to
$\p_{\eps}.$
We will now show that
from the convexity of $\Phi$ and $\Phi^{\prime}$
and Proposition 1 it follows 
\begin{equation}
\p_{\eps}(
\Phi(\eps)\geq M +
\|\Phi\|_{L}\sqrt{t} 
)\leq 2^{\alpha}
\exp\Bigl\{-\frac{\alpha}{\alpha+1}t\Bigr\},
\label{star}
\end{equation}
and
\begin{equation}
\p_{\eps}(
\Phi^{\prime}(\eps)\leq M -
\|\Phi^{\prime}\|_{L}\sqrt{t} 
)\leq 2^{\alpha}
\exp\Bigl\{-\frac{\alpha}{\alpha+1}t\Bigr\}.
\label{2star}
\end{equation}
Let us recall how this is usually done (see Ledoux and Talagrand (1991)).
If we consider the set 
${\cal A}=\{\eps : \Phi(\eps)\leq M\},$ then 
$\P({\cal A})\geq 1/2$ and by convexity of $\Phi,$
$\mbox{conv}{\cal A}= {\cal A}.$
This, together with the Lipschitz condition, implies that
$$
\{f_c^2({\cal A}, \eps)\leq t\}\subseteq
\{\Phi(\eps)\leq M +\|\Phi\|_{L}\sqrt{t}\}.
$$
Thus, the right tail (\ref{star}) follows from Proposition 1.
Similarly, if we consider the set 
$$
{\cal B}=\{\eps : \Phi^{\prime}(\eps)\leq M
- \|\Phi^{\prime}\|_{L}\sqrt{t} \},
$$ 
then 
$$
\{f_c^2({\cal B}, \eps)\leq t\}\subseteq
\{\Phi^{\prime}(\eps)\leq M \}.
$$
By Proposition 1,
$$
\frac{1}{2}\leq 
\p( f_c^2({\cal B}, \eps)\geq t)\leq \frac{1}{\p({\cal B})^{\alpha}}
\exp\Bigl\{-\frac{\alpha}{\alpha+1}t\Bigr\}.
$$
We can rewrite this as
$$
\p({\cal B})\leq 2^{\beta}
\exp\Bigl\{-\frac{\beta}{\beta+1}t\Bigr\},
$$
where $\beta=1/\alpha.$ But since $\alpha$ is arbitrary, this proves
the lower tail (\ref{2star}),
which completes the proof of the theorem.
\qed

This result is an intermediate step in obtaining the
concentration inequality for $Z(x)$ in its final form,
since $V$ still depends on $x.$
Notice that here we did not assume any boundedness of
$f\in\F,$ and the result is of somewhat similar nature
as the self-normalization phenomenon in the one-dimensional case
(see Gin\'e et. al. (1997), or Shao(1997)).
Under the additional assumption that $f\in\F$ are uniformly bounded
one can proceed by controlling the deviation of
$V$ (or $W$) from its expectation,
which is done in a usual way, either via control by
two points as in Talagrand (1996) plus some truncation argument, 
or via a sharp concentration inequality of Boucheron et. al. (2000).

Let us assume now that 
$$
\forall f\in\F\ \forall x\in\Omega,\,\,\,
-\frac{1}{2}\leq f(x)\leq \frac{1}{2}.
$$
If we introduce
$V_i=\e_{y}\sup_{f\in\F}\sum_{j\not = i}(f(x_j)-f(y_j))^2$
then, it is easy to see that
$$
0\leq V-V_i\leq 1\,\,\, \mbox{ and }\,\,\,
\sum_{i=1}^n (V-V_i)\leq V.
$$
Under these conditions, Theorem 6 in Boucheron et. al. (2000) states
that for all $t\geq 0,$
\begin{equation}
\p\Bigl(V\geq \e V +t\Bigr)\leq \exp\Bigl\{-\e V
h\Bigl(\frac{t}{\e V}\Bigr) \Bigr\},
\label{Pois}
\end{equation}
where $h(x)=(1+x)\log(1+x) -x.$
Since 
$h(x)\geq x^2/(2+2x/3),$
(\ref{Pois}) implies Bernstein's inequality
$$
\p\Bigl(V\geq \e V +t\Bigr)\leq \exp\Bigl\{
-\frac{t^2}{2\e V + 2t/3} \Bigr\},
$$
which can be equivalently written as
$$
\p\Bigl(
V\geq \e V+ \frac{1}{3}
(18\e V t +t^2)^{1/2}+\frac{t}{3}
\Bigr)\leq e^{-t}.
$$
More generally, if $-b\leq f(x)\leq b,$
then
$$
\p\Bigl(
V\geq \e V+ \frac{2b}{3}
(18\e V t +4b^2 t^2)^{1/2}+\frac{4b^2t}{3}
\Bigr)\leq e^{-t}.
$$
Combining this with Theorem 1 we get the following corollary.
\begin{corollary}
If $-b\leq f(x)\leq b$ then for all $t\geq \log 2,$
\begin{equation}
\p\Bigl(
|Z-\e Z|\geq
2\Bigl(
t\Bigl(\e V+ \frac{2b}{3}
(18\e V t +4b^2 t^2)^{1/2}+\frac{4b^2t}{3}
\Bigr)
\Bigr)^{1/2}
\Bigr)
\leq 
4e^{1-(\sqrt{t}-\sqrt{\log 2})^2} +e^{-t}.
\label{tails}
\end{equation}
\end{corollary}
\qed

It is clear, that in the range of parameters
$1\ll t\ll \e V/b^2,$ the bound of the Corollary will be dominated
by the term $\sim 2\sqrt{\e V t}.$ 
For this range, it improves upon the control
of the lower tail given by Theorem 12 in  Massart (2000), which states
\begin{equation}
\p\Bigl(
Z\leq \e Z -2\sqrt{1.35 \e V t} -3.5 b t
\Bigr)\leq e^{-t}.
\label{tail}
\end{equation}
Actually, one can check that 
$$
2\Bigl(
t\Bigl(\e V+ \frac{2b}{3}
(18\e V t +4b^2 t^2)^{1/2}+\frac{4b^2t}{3}
\Bigr)
\Bigr)^{1/2}
\leq 
2\sqrt{1.35 \e V t} +3.5 b t
$$
for all parameters $b,\e V, t.$ 
Unfortunately, (\ref{tails}) and (\ref{tail})
are not comparable in all range of parameters, mainly,
because of the term $\exp\{-(\sqrt{t}-\sqrt{\log 2})^2\}.$

Finally, for more results in this 

\vspace{0.5cm}

\section{Vapnik-Chervonenkis type inequalities.}

In this section we are trying to control the functional $Q_n f$ 
uniformly over the class $\F,$ where
$$
Q_n f=Pf-P_n f \,\,\,\mbox{ or }\,\,\, Q_n f = P_n f -P f
$$
and
$$
Pf=\int f(x) dP(x),\,\,\,
P_n f = \frac{1}{n}\sum_{i=1}^{n} f(x_i).
$$
The difference from the previous section is that now 
the bounds on $Q_n f$ will depend on $f$ and will reflect
that the function $f$ with a smaller variance should have
a tighter bound. The results of this section
are in a spirit of Vapnik and Chervonenkis (1968) and 
Panchenko (2002).

Corresponding to $Q_n f$, let us introduce
$$
S_n f=\frac{1}{n}\sum_{i=1}^n (f(y_i)-f(x_i))
\,\,\,\mbox{ or }\,\,\, 
S_n f=\frac{1}{n}\sum_{i=1}^n (f(x_i)-f(y_i)).
$$
Finally, we define
$$
R_n f=\frac{1}{n}\sum_{i=1}^n \eps_i (f(y_i)-f(x_i)),
$$
$$
W f=W(f,x,y)=\frac{4}{n}\sum_{i=1}^n (f(y_i)-f(x_i))^2,\,\,\,\,
V f=V(f, x)=\e_{y}W(f, x, y).
$$
As one of the consequences of our approach we will give a uniform
control of $Q_n f/ (V f)^{1/2}$ for VC-subgraph classes of functions.
The original result of Vapnik and Chervonenkis \cite{VaCher}
provided a uniform control for $Q_n f/ (P f)^{1/2}$ for
VC-classes of functions taking values $f\in\{0,1\}$ 
(and a simple generalization for VC-major classes taking values
in $[0,1]$). The fact that we can substitute $P f$ by $V f$  
gives a new way to control $Q_n f.$ 

Let us introduce a function $\Phi(f,x,y)$ which is invariant over all
permutations of $(x,y)$ that switch only the same coordinates
of $x$ and $y.$
Assume that for some fixed $\beta\in(0,1)$ and for any fixed $(x,y)$
we have
\begin{equation}
\p_{\eps}\Bigl(
\sup_{f\in\F}(R_n f -\Phi(f,x,y))> 0
\Bigr) < 1-\beta.
\label{Phi}
\end{equation}
Then the following theorem holds.
\begin{theorem}
Assume that (\ref{Phi}) holds. Then for any $t\geq \log\beta^{-1},$
$$
\p\Bigl(
\exists f\in\F \,\,
Q_n f\geq \e_{y} \Phi(f,x,y)+\sqrt{\frac{Vt}{n}}
\Bigr)\leq
\exp(1-(\sqrt{t}-\sqrt{\log\beta^{-1}})^2).
$$
\end{theorem}
{\bf Proof.} 
We will first prove that for any $\alpha\geq 0$ 
the statement of the theorem holds with the right hand
side substituted by $\beta^{-\alpha}\exp(1-\alpha t/(\alpha+1)).$
The result will follow by optimization over $\alpha.$
First of all, by Corollary 1 it is enough to prove 
that
$$
\p\Bigl(
\exists f \,\,
S_n f\geq \Phi(f,x,y)+\sqrt{\frac{Wt}{n}}
\Bigr)\leq
\frac{1}{\beta^{\alpha}}
\exp\Bigl(-\frac{\alpha}{\alpha +1}t\Bigr).
$$
Since
$\Phi(f,x,y)$ is invariant under permutations of $x_i$ and $y_i.$
we can write,
\begin{eqnarray}
&&
\p\Bigl(
\exists f \,\,
S_n f\geq \Phi(f,x,y)+\sqrt{\frac{Wt}{n}}
\Bigr)=
\p\Bigl(
\exists f \,\,
R_n f\geq  \Phi(f,x,y)+\sqrt{\frac{Wt}{n}}
\Bigr)
\nonumber
\\
&&
=\e\p_{\eps}\Bigl(
\exists f \,\,
R_n f\geq \Phi(f,x,y)+\sqrt{\frac{Wt}{n}}
\Bigr).
\label{symright}
\end{eqnarray}
For a fixed $(x,y)$ consider a set
$$
{\cal A}=\{\eps : \sup_{f\in\F}(R_n f - \Phi(f,x,y))\leq 0\}.
$$
By condition (\ref{Phi}), $\p_{\eps}({\cal A})\geq\beta.$
If we denote 
${\cal A}_t=\{\eps : f_c^2({\cal A},\eps)\leq t\}$
then (\ref{T1}) implies that 
$$
\p_{\eps}({\cal A}_t)\geq 1-\beta^{-\alpha}
\exp\Bigl(-\frac{\alpha}{\alpha+1}t\Bigr).
$$
Let us take $\eps\in{\cal A}_t$ and $\eps'\in{\cal A}.$
The definition of ${\cal A}$ implies that for any $f\in\F$
$$
\frac{1}{n}\sum_{i=1}^{n}\eps_i'(f(y_i)-f(x_i))\leq \Phi(f,x,y),
$$
and, therefore,
\begin{eqnarray*}
&&
\frac{1}{n}\sum_{i=1}^{n}\eps_i(f(y_i)-f(x_i)) - \Phi(f,x,y)
\leq 
\frac{1}{n}\sum_{i=1}^{n}(\eps_i-\eps_i')(f(y_i)-f(x_i))
\\
&&
\leq
\frac{2}{n}\sum_{i=1}^{n}|f(y_i)-f(x_i)|I(\eps_i'\not =\eps_i).
\end{eqnarray*}
But since $\eps\in{\cal A}_t,$ (\ref{T2}) implies that
one can choose $\eps'\in{\cal A}$ so that
$$
\frac{2}{n}\sum_{i=1}^{n}|f(y_i)-f(x_i)|I(\eps_i'\not =\eps_i)
\leq
\Bigl(
t\frac{4}{n^2}\sum_{i=1}^{n}(f(y_i)-f(x_i))^2
\Bigr)^{1/2}=
\Bigl(
\frac{Wt}{n}
\Bigr)^{1/2}.
$$
This proves the theorem.
\qed

Let us consider a special case of $\Phi(f,x,y),$ 
which satisfies condition (\ref{Phi}). 
Let us note here that application of Talagrand's concentration 
inequality for two point space
as it was implemented in Theorem 2 is not
crucial for the examples of this section. 
It is well known fact that the chaining technique that
we will only use here to bound the $(1-\beta)$-quantile implies
tail estimates as well. But it is hard to argue with the fact that
the application of Talagrand's inequality even for these examples
is more elegant as it immediately provides the tail estimates once the
bound for the quantile is obtained.

We will assume from now on that $0\equiv f\in \F.$
Let $d$ be a metric on $\F.$
Given $u>0$ we say that a subset ${\cal F}'\subset {\cal F}$ 
is $u-$separated if for any $f\not =g\in {\cal F}'$ we have
$ d(f,g)> u.$ Let a {\it packing number} 
$D({\cal F},u,d)$ be the maximal cardinality of 
a $u-$separated set. 

We define
$$
\Phi(f,x,y)=Kn^{-1/2}\int_{0}^{\sqrt{W}/2}(\log D(\F,u,d_{x,y}))^{1/2}du,
$$
where 
$$
d_{x,y}(f,g)=\Bigl(
\frac{1}{n}\sum_{i=1}^{n}(f(y_i)-f(x_i)-g(y_i)+g(x_i))^2
\Bigr)^{1/2}
$$
and $K=K(\beta)$ depends only on $\beta.$ 
For example, if $K(\beta)=8(p+2)^{1/2},$ where $p$ is such that
$\sum_{j=2}^{\infty}j^{-p}<1-\beta,$ then the following theorem holds.
\begin{theorem}
If $K(\beta)$ is defined as above then (\ref{Phi}) holds.
\end{theorem}
{\bf Proof.} 
The proof is based on standard chaining technique.
Let us fix $(x,y).$ Define
$$
F=\{(f(y_1)-f(x_1),\ldots,f(y_n)-f(x_n)) : f\in\F\}
$$
and
$$
d(f,g)=\Bigl(
\frac{1}{n}\sum_{i=1}^{n}(f_i-g_i)^2\Bigr)^{1/2},\,\,\,\,
f,g\in F.
$$
Then, if 
$$
\Phi(f)=K(\beta)n^{-1/2}\int_{0}^{d(f,0)}(\log D(F,u,d))^{1/2}du,
$$
we need to prove that
$$
\p_{\eps}\Bigl(\sup_{f\in F}\Bigl(\frac{1}{n}
\sum_{i=1}^{n}\eps_i f_i - \Phi(f)
\Bigr)
> 0)< 1-\beta.
$$
Let $j_0$ be defined as
$$
j_0=\inf\{j : D(F,2^{-j},d)\geq 2\}.
$$
Consider an increasing sequence of sets
$$
\{0\}=F_{-\infty}=\ldots=F_{j_0-1}
\subseteq F_{j_0}\subseteq F_{j_0+1}\subseteq \ldots
$$
such that for any $g\not =h\in F_j,$ $d(g,h)>2^{-j}$
and for all $f\in F$ there exists $g\in F_j$ such that
$d(f,g)\leq 2^{-j}.$ The cardinality of $F_j$ can be bounded by
$$
|F_j|\leq D(F,2^{-j},d).
$$
For simplicity of notations we will write $D(u):= D(F,u,d).$
If $D(2^{-j})=D(2^{-j-1})$ then in the construction of the sequence
$(F_j)$ we will set $F_j$ equal to $F_{j+1}.$
We will now define the sequence of projections 
$\pi_j:F\to F_j,\,j\geq 0$ in the following way.
If $f\in F$ is such that $d(f,0)\in(2^{-j-1},2^{-j}]$
then set $\pi_0(f)=\ldots=\pi_j(f)=0$ and for $k\geq j+1$ 
choose $\pi_k(f)\in F_k$ such that $d(f,\pi_k(f))\leq 2^{-k}.$
In the case when $F_k=F_{k+1}$ we will choose $\pi_{k}(f)=\pi_{k+1}(f).$
This construction implies that $d(\pi_{k-1}(f),\pi_{k}(f))\leq 2^{-k+2}.$
Let us introduce a sequence of sets
$$
\Delta_j=\{g-h : g\in F_j, h\in F_{j-1}, d(g,h)\leq 2^{-j+2}\},\,\,\,\,
j\geq j_0,
$$
and let $\Delta_j = \{0\}$ if $D(2^{-j})=D(2^{-j+1}).$
The cardinality of $\Delta_j$ does not exceed 
$$
|\Delta_j|\leq |F_j|^2\leq D(2^{-j})^2.
$$
By construction any $f\in F$ can be represented 
as a sum of elements from $\Delta_j$
$$
f=\sum_{j\geq j_0}(\pi_j(f)-\pi_{j-1}(f)),\,\,\,\,\,\,
\pi_{j}(f)-\pi_{j-1}(f) \in \Delta_j.
$$
Let
$$
I_j=n^{-1/2}\int\limits_{2^{-j-1}}^{2^{-j}}(\log D(u))^{1/2} du
$$
and define the event
$$
A=
\bigcup_{j=j_0}^{\infty}
\Bigl\{
\sup_{f\in\Delta_j}\frac{1}{n}\sum_{i=1}^{n}\varepsilon_i f_{i}
\geq K I_j
\Bigr\}.
$$
On the complement $A^c$ of the event $A$ we have for any 
$f\in F$ such that $d(f,0)\in(2^{-j-1},2^{-j}]$
\begin{eqnarray*}
&&
\frac{1}{n}\sum_{i=1}^{n}\varepsilon_i f_{i}=
\frac{1}{n}\sum_{k\geq j+1}\sum_{i=1}^{n}\varepsilon_i 
(\pi_k(f)-\pi_{k-1}(f))_{i}
\leq
\sum_{k\geq j+1} K I_k 
\\
&&
\leq
Kn^{-1/2}\int\limits_{0}^{2^{-j-1}}(\log D(u))^{1/2} du\leq
Kn^{-1/2}\int\limits_{0}^{d(f,0)}(\log D(u))^{1/2} du.
\end{eqnarray*}
It remains to prove that for some constant $K(\beta),$
$P(A)<1-\beta.$ Indeed,
\begin{eqnarray*}
&&
P(A)\leq \sum_{j=j_0}^{\infty}
P\Bigl(
\sup_{f\in\Delta_j} \frac{1}{n}\sum_{i=1}^{n}\varepsilon_i f_{i}
\geq K I_j
\Bigr) 
\\
&&
\leq
\sum_{j=j_0}^{\infty}
|\Delta_j|\exp \Bigl\{
-\frac{n K^2 I_j^2}{2^{-2j+4}}
\Bigr\}
I\bigl(D(2^{-j})>D(2^{-j+1})\bigr)
\\
&&
\leq
\sum_{j=j_0}^{\infty}
\exp \Bigl\{
2\log D(2^{-j})
-\frac{n K^2 I_j^2}{2^{-2j+4}}
\Bigr\}
I\bigl(D(2^{-j})>D(2^{-j+1})\bigr),
\end{eqnarray*}
since for $f\in\Delta_j$
$n^{-1}\sum_{i=1}^{n}f_i^2\leq 2^{-2j+4}.$
The fact that $D(u)$ is decreasing implies
$$
\frac{n^{1/2}I_j}{ 2^{-(j+1)}}\geq (\log D(2^{-j}))^{1/2}
$$
and, therefore,
\begin{eqnarray*}
&&
P(A)\leq \sum_{j=j_0}^{\infty}
\exp \{
-\log D(2^{-j})(K^2 2^{-6}-2)
\}
I\bigl(D(2^{-j})>D(2^{-j+1})\bigr)
\\
&&
\leq
\sum_{j=j_0}^{\infty}
\frac{1}{D(2^{-j})^{p}}
I\bigl(D(2^{-j})>D(2^{-j+1})\bigr)
\leq
\sum_{j=2}^{\infty}
\frac{1}{j^{p}}<1-\beta,
\end{eqnarray*}
for $p=K(\beta)^2 2^{-6} -2$ big enough. We used the fact
that $D(2^{-j_0})\geq 2.$
\qed

{\bf \hspace{-1cm} Example} (Uniform entropy conditions).
Let us introduce a uniform packing numbers
$D(\F,u)$ as any function such that 
$$
\sup_{Q}D({\cal F},u,L_2(Q))\leq D({\cal F},u)
$$
where the supremum is taken over all discrete probability measures.
One can easily check that
$$
\Bigl(\frac{1}{n}\sum_{i=1}^{n}
(f(x_i)-f(y_i)-g(x_i)+g(y_i))^2\Bigr)^{1/2}\leq
2
\Bigl(\frac{1}{2n}\sum_{i=1}^{n}\bigl(
(f(x_i)-g(x_i))^2+(f(y_i)-g(y_i))^2\bigr)
\Bigr)^{1/2}
$$
and, therefore, in the case 
when the packing numbers are bounded uniformly we get,
$$
D(\F,u,d_{x,y})\leq D(\F,u/2).
$$
Hence,
\begin{eqnarray*}
&&
\e_{y}\Phi(f,x,y)\leq
K(\beta)n^{-1/2}\e_{y}\int_{0}^{\sqrt{W}/2}(\log D(\F,u/2))^{1/2}du
\\
&&
\leq 2K(\beta)n^{-1/2}\int_{0}^{\sqrt{V}/4}(\log D(\F,u))^{1/2}du.
\end{eqnarray*}
\begin{corollary}
For any $t\geq \log\beta^{-1},$
$$
\p\Bigl(
\exists f\in\F \,\,
Q_n f \geq
\frac{2K(\beta)}{n^{1/2}}\int_{0}^{\sqrt{V}/4}(\log D(\F,u))^{1/2}du
+\sqrt{\frac{Vt}{n}}
\Bigr)\leq
\exp(1-(\sqrt{t}-\sqrt{\log\beta^{-1}})^2).
$$
\end{corollary}
\qed

In the case of VC-subgraph classes with VC dimension $d$
(for definition, see van der Vaart and Wellner (1996)), the result of 
\cite{Haussler} gives
$$
D(\F,u)\leq e(d+1)\Bigl(\frac{2e}{u^2}\Bigr)^d,
$$
and, therefore, the following corollary.
\begin{corollary}
(Normalization by variance).
There exists $K$ that depends only on $\beta$ such that
for any $t\geq \log\beta^{-1},$
$$
\p\Bigl(
\exists f\in\F \,\,
\frac{Q_n f}{\sqrt{V}} \geq
K\sqrt{\frac{d\log n}{n}}
+\sqrt{\frac{t}{n}}
\Bigr)\leq
\exp(1-(\sqrt{t}-\sqrt{\log
\beta^{-1}})^2).
$$
\end{corollary}
\qed

Let us rewrite $V$ as
$$
V=V(x)=4(\mbox{Var}f+\mbox{Var}_n f +(Pf-P_n f)^2)
=4(\mbox{Var}f+\mbox{Var}_n f +(Q_n f)^2) ,
$$
where
$$
\mbox{Var}_n f=\frac{1}{n}\sum_{i=1}^{n}(P_n f - f(x_i))^2
$$
is a sample variance.
If we denote
$$
U=K\sqrt{\frac{d\log n}{n}}+\sqrt{\frac{t}{n}},
$$
then one can solve the inequality of Corollary 4 for $Q_n f$ to get
$$
\p\Bigl(
\exists f\in\F \,\,
|Q_n f| \geq
2U\Bigl(\frac{\mbox{Var}f+\mbox{Var}_n f}{1-4U^2}\Bigr)^{1/2}
\Bigr)\leq
2\exp(1-(\sqrt{t}-\sqrt{\log\beta^{-1}})^2).
$$

Let us compare this to an ``optimistic'' inequality of 
Vapnik and Chernonenkis \cite{Vap}, which states that if 
$\F=\{f:\Omega\to\{0,1\}\}$ is a VC-class of indicator
functions with VC dimension $d$, then with probability
at least $1-e^{-t/4},$ for all $f\in \F$
$$
\frac{1}{n(Pf)^{1/2}}\sum_{i=1}^{n}(Pf - f(x_i))\leq
2\Bigl(
\frac{d }{n}\log \frac{2en}{d} + 
\frac{t}{n}
\Bigr)^{1/2}.
$$
Compared to the
inequality of Vapnik and Chervonenkis our inequality
controls the deviation of $P_n f$ from $P f$ in both directions,
no assumptions are made on the boundedness of functions $f\in \F,$
and the deviation is controled by the mixture of variance and sample
variance rather than by expectation $P f,$ which can be considered
as a significant improvement.

{\bf Example} (The case of one function).
When $\F$ consists of one function $f$ we will simply write
$f(X)=\xi.$ Let us take $\beta=1/2$ and let
$$
\Phi(\xi)=\e_{\xi'}
M_{\eps}\Bigl(\frac{1}{n}\sum_{i=1}^n \eps_i(\xi_i-\xi_i')\Bigr)=0.
$$
Obviously, with this choice of $\beta$ and $\Phi$ 
condition (\ref{Phi}) holds and Theorem 2 implies
$$
\p\Bigl(|\bar{\xi}-\e\xi|\geq
2\Bigl(\frac{(\mbox{Var}\xi+\mbox{Var}_n \xi +(\e\xi-\bar{\xi})^2)t}{n}
\Bigr)^{1/2}
\Bigr)\leq 2\exp\Bigl(1-(\sqrt{t}-\sqrt{\log 2})^2\Bigr).
$$
Solving the inequality for $|\bar{\xi}-\e\xi|$ we get
\begin{equation}
\p\Bigl(|\bar{\xi}-\e\xi|\geq
2\Bigl(\frac{(\mbox{Var}\xi+\mbox{Var}_n \xi)t}{n-4t}\Bigr)^{1/2}
\Bigr)\leq 2\exp\Bigl(1-(\sqrt{t}-\sqrt{\log 2})^2\Bigr).
\label{one}
\end{equation}
One should compare this to Bernstein type inequalities. 
First of all, we don't assume any moment conditions other than
the existance of variance of $\xi.$ 
Second, in Bernstein's inequality
$$
|\bar{\xi}-\e\xi|\lessim \Bigl(\frac{t\mbox{Var}\ \xi }{n}\Bigr)^{1/2}
\mbox{ for } t\leq n \mbox{Var}\ \xi,
$$
$$
|\bar{\xi}-\e\xi|\lessim \frac{t}{n}
\mbox{ for } t\geq n \mbox{Var}\ \xi,
$$
whereas (\ref{one}) gives
$$
|\bar{\xi}-\e\xi|\leq 2\Bigl(\frac{2(\mbox{Var}\xi+\mbox{Var}_n\xi)t}
{n}\Bigr)^{1/2}
\mbox{ for } t\leq n/8.
$$
This, basically, means that the deviation of the average
$\bar{\xi}$ from the expectation $\e\xi$ can be large only
when the sample variance is large.
\vspace{0.5cm}

{\bf Acknowledgment.} We want to thank Michel Talagrand for
some valuable comments and suggestions.

\vskip 1mm

\hfill\break
Department of Mathematics and Statistics\hfill\break
The University of New Mexico\hfill\break
Albuquerque, NM 87131--1141\hfill\break
e-mail: panchenk@math.unm.edu\hfill\break


\begin{thebibliography}{99}

\bibitem{Ma1} Boucheron, S., Lugosi, G., Massart, P. (2000)
A sharp concentration inequality with applications. 
\it Random Structures Algorithms \rm {\bf 16}  277 - 292. 


\bibitem{Dem} Dembo, A. (1997)
Information inequalities and concentration
of measure. \it Ann. Probab. \rm {\bf 25} 527 - 539.

\bibitem{Du} Dudley, R.M. (1999) Uniform Central Limit Theorems.
Cambridge University Press.

\bibitem{Gine} Gin\'e, E., G\"otze, F., Mason, D. (1997)
When is the Student $t$-statistics asymptotically standard normal?
{\it Ann. Probab.} {\bf 26} 1514 - 1431.

\bibitem{Haussler}
Haussler, D. (1995)
Sphere packing numbers for subsets
of the boolean $n-$cube with bounded Vapnik-Chervonenkis dimension. 
\it J. Combin. Theory Ser. A \rm {\bf 69} 217 - 232.

\bibitem{Le} Ledoux, M. (1996)
On Talagrand's deviation inequalities for product measures.
{\em ESAIM: Probab. Statist.} {\bf 1} 63 - 87. 

\bibitem{LT} Ledoux, M. and Talagrand, M. (1991)
Probability in Banach Spaces.
Springer-Verlag, New York.

\bibitem{Ma} Massart, P. (2000)
About the constants in Talagrand's concentration inequalities
for empirical processes. {\em Ann. Probab.} {\bf 28} 863 - 885.

\bibitem{Pa} Panchenko, D. (2001)
A note on Talagrand's concentration inequality.
\it Elect. Comm. in Probab. \rm {\bf 6} 55 - 65. 

\bibitem{Pa1} Panchenko, D. (2002)
Some extensions of an inequality of Vapnik and Chervonenkis. 
\it Elect. Comm. in Probab. \rm {\bf 7} 

\bibitem{Rio} Rio E. (2000) 
In\'egalit\'es exponentielles pour
les processus empiriques.
\it C.R. Acad. Sci. Paris, \rm t.330, S\'erie I 597-600.

\bibitem{Rio2} Rio E. (2001)
In\'egalit\'es de concentration pour
les processus empiriques de classes de parties. 
\it Probab. Theory Relat. Fields \rm {\bf 119} 163-175.

\bibitem{Shao} Shao, Q.-M. (1997)
Self-normalized large deviations. 
{\em Ann. Probab.} {\bf 25} 285 - 329.

\bibitem{Ta1} Talagrand, M. (1995)
Concentration of measure and
isoperimetric inequalities in product spaces.  
\it Publications Math\'ematiques de l'I.H.E.S. \rm {\bf 81} 73-205.

\bibitem{Ta} Talagrand, M. (1996)
New concentration inequalities in product spaces. {\em Invent.
Math.} {\bf 126} 505-563.

\bibitem{Well} van der Vaart, A., Wellner, J. (1996)
Weak Convergence and Empirical Processes: With Applications to Statistics.
John Wiley \& Sons, New York. 

\bibitem{VaCher} Vapnik, V.N., Chervonenkis, A.Ya. (1968)
On the uniform convergence of relative frequencies of event to
their probabilities.
\it Soviet Math. Dokl.\rm {\bf 9} 915 - 918.

\bibitem{Vap} Vapnik, V.N. (1998)
Statistical Learning Theory. Wiley, New York.

\end{thebibliography}
\end{document}